\def\version{July 29, 2008}

\documentclass[12pt]{article}
\usepackage[psamsfonts]{amsfonts}
\usepackage{amsmath,amssymb,amsthm}
\newtheorem{theorem}{Theorem}[section]
\newtheorem{prop}[theorem]{Proposition}
\newtheorem{lemma}[theorem]{Lemma}

\newtheorem{defn}[theorem]{Definition}
\newtheorem{rem}[theorem]{Remark}
\newtheorem{ass}[theorem] {Assumption}
\makeatletter
\@addtoreset{equation}{section}

\makeatother

\def\eq#1\en{\begin{equation}#1\end{equation}}  
\newcommand{\nnb}	{\nonumber \\} 
\def\eqalign#1\enalign{
	\begin{align}#1\end{align}}
\newcommand{\lbeq}[1]  {\label{e:#1}}
\newcommand{\refeq}[1] {\eqref{e:#1}}    

\newcommand{\Ccal}   {\mathcal{C}} 
 
\newcommand{\Ecal}   {\mathcal{E}} 
\newcommand{\Fcal}   {\mathcal{F}}

\newcommand{\bE}{{\mathbb E}}
\newcommand{\bN}{{\mathbb N}}
\newcommand{\bP}{{\mathbb P}}
\newcommand{\bR}{{\mathbb R}}
\newcommand{\bZ}{{\mathbb Z}}

\newcommand{\IIC}{{\hbox{\footnotesize\rm IIC}}}
\def\Reff{R_{\rm eff}}  
\def\tReff{{\tilde R}_{\rm eff}}
\def\eps{\varepsilon}
\def\lam{\lambda}
\def\Gam{\Gamma}
\newcommand{\al}{\alpha}

\renewcommand{\to}      {\rightarrow}

\title  {Heat kernel estimates for strongly recurrent
random walk on random media} 

\author{Takashi Kumagai
\footnote{Department of Mathematics, 
Kyoto University, Kyoto 606-8502, Japan.
{\tt kumagai@math.kyoto-u.ac.jp}}
\thanks{Research partially supported by the
Grant-in-Aid for Scientific Research (B) 18340027.}
and Jun Misumi
\footnote{Graduate School of Mathematical Sciences,
The University of Tokyo, Komaba, Tokyo 153-8914, Japan.
{\tt misumi@ms.u-tokyo.ac.jp}}
\thanks{Research partially supported 
by the 21 century COE program
at Graduate School of Mathematical Sciences, the University of Tokyo.}
}
\begin{document}
\date\version
\maketitle
\begin{abstract}
We establish general estimates
for simple random walk on an arbitrary infinite
random graph, assuming suitable bounds on volume and effective
resistance for the graph. These are generalizations of the results in 
\cite[Section 1,2]{BJKS}, 
and in particular, imply the spectral dimension of the random graph. 
We will also give an application of the results to random walk on  
a long range percolation cluster. \\
{\it Key Words}: Random walk \,-\, Random media \,-\, 
Heat kernel estimates \,-\, Spectral dimension \,-\, Long range percolation\\
{\it Running Head}: Heat kernel estimates on random media
\end{abstract}

\section{Introduction and Main results}
\subsection{Introduction}
Recently, there are intensive study for detailed properties of 
random walk on a percolation cluster. For a random walk on a 
supercritical percolation cluster on $\bZ^d$\,($d\ge 2$), detailed Gaussian
heat kernel estimates and quenched invariance principle have
been obtained (\cite{Barl04,BB06,MP06,SS04}).
This means, such a random walk behaves in a diffusive fashion similar 
to a random walk on $\bZ^d$. On the other hand,
it is generally believed that random walk on a large 
critical cluster behaves subdiffusively (see \cite{HBA00} and the references therein).  
Critical percolation clusters are believed (and for some cases proved) 
to be finite. So, it is natural to consider random walk
on an incipient infinite cluster (\IIC), namely a critical percolation
cluster conditioned to be infinite. Random walk on \IIC{}s 
has been 
proved to be subdiffusive on $\bZ^2$ 
(\cite{Kest86}), on trees (\cite{Kest86a,BK06}), 
and for the spread-out oriented percolation
on $\bZ^d\times \bZ_+$ in dimension $d>6$ (\cite{BJKS}). 

In order to study detailed properties of the random walk, it is 
nice and useful if one can compute the long time behaviour of the transition
density (heat kernel). Let $p_n(x,y)$ be its transition density (see 
\refeq{TAdef11} for a precise definition) of a random walk on an infinite 
(random) graph $G$. Define the {\sl spectral dimension} of  $G$ by
$$ d_s(G) =-2 \lim_{n \to \infty} \frac{\log p_{2n}(x,x)}{\log n},$$
(if this limit exists). Let ${\cal C}_d$ be the \IIC\, for percolation
cluster on $\bZ^d$. 
Alexander and Orbach \cite{AO82} conjectured that, for any $d\ge 2$, 
$d_s({\cal C}_d)=4/3$. (To be precise, the original 
Alexander-Orbach conjecture was that the left side hand of 
(\ref{e:oaocpnj}) is equal to $2/3$ for all dimensions on ${\cal C}_d$.) 
While it is now believed that this conjecture is unlikely
to be true for small $d$ (\cite{HBA00}), it is proved that the conjecture is true
on trees (\cite{Kest86a,BK06}), and for the spread-out oriented percolation
on $\bZ^d\times \bZ_+$, $d>6$ (\cite{BJKS}). 
One very interesting open problem in this direction is to establish estimates
of $d_s({\cal C}_d)$ (or the corresponding \IIC\, for the oriented percolation)
for $d$ small to disprove/prove the Alexander-Orbach conjecture. 

In \cite{BJKS}, general estimates are given for simple random walk on an 
arbitrary infinite random graph $G$, assuming suitable bounds on volume and 
effective resistance for the random graph.  In particular, $d_s(G)=4/3$. 
The main purpose of this paper is to extend this  
general estimates to the framework of 
strongly recurrent random walk on the random graph $G$. Here 
`strongly recurrent' simply means $d_s(G)<2$. (See \cite{BCK05} Section 1.1,
for more precise meaning of `strongly recurrent'.) Roughly saying our main
results can be expressed as follows; if the volume of the ball of radius $R$ 
on the random graph is of order $R^D$ with high probability and the resistance
between the center and the outside of the ball is of order $R^\alpha$ with high  probability (precise statement is given 
in Assumption \ref{ass-rwre}), then one can establish 
both quenched (i.e. almost sure with respect to the randomness of the graph) 
and annealed (i.e. averaged over the randomness of the graph) estimates for the exit time from the ball, on-diagonal heat kernel, the mean displacement, etc.
(Propositions~\ref{ptight}--\ref{pmeans} and 
Theorem~\ref{thm-rwre}). 
In particular,
\begin{equation}\label{e:aldds}
d_s(G)=\frac{2D}{D+\alpha}.\end{equation}
Note that the estimates given in Section 1--2 of 
\cite{BJKS} treat the case $D=2,\alpha=1$.   
Our results are general enough to allow logarithmic corrections
for the volume and the resistance estimates. In fact, the volume growth  
$v(R)$ and the resistance growth $r(R)$ could be any functions satisfying 
(\ref{vrdoubl}). 

Unfortunately, so far we could not establish any new results for random walk on the 
\IIC\, for low dimensions. Though our results give `simplest' conditions
on volume and resistance growth (Assumption \ref{ass-rwre}) 
to obtain the spectral dimension $d_s(G)$, significant work is required to prove 
them for such models, as was done in Section 3--5 of \cite{BJKS} 
for the spread-out oriented percolation on $\bZ^d\times \bZ_+$, $d>6$. 
Instead, we apply our results to the long range percolation, which is another
important random graph. We consider the following case of the long range percolation. 
On $\bZ$, each pair of points $x,y\in\bZ$, $|x-y|\ge 2$ is
connected by an unoriented bond with probability $\beta|x-y|^{-s}$ for some 
$\beta>0$ and $s>2$, independently of each others. Each pair of nearest
points is connected with probability $1$. On this random graph, we can check
the  bounds on volume and effective resistance required for the general estimates
and obtain, for example, (\ref{e:aldds}) with $D=\alpha=1$ (Theorem \ref{lr-result}).   
In fact, it is quite likely that for $s>2$, the transition density for simple 
random walk on the long range percolation is   
Gaussian-type, so there are other ways to establish (\ref{e:aldds}). 
However, we can further 
observe an interesting discontinuity of the spectral dimension 
at $s=2$ (Remark \ref{lr-rem}(1)). 

Finally, we remark that after the first draft of this paper was submitted,  
some rigorous results were obtained concerning the Alexander-Orbach conjecture. 
In \cite{KozNach}, it is proved that the conjecture is true when the two-point 
function behaves nicely; in particular the conjecture is true  
for simple random walk on the \IIC{} of critical percolation on $\bZ^d$ when 
$d$ is large enough or when $d>6$ and the lattice is sufficiently spread out. 
In \cite{CroKum}, it is proved that the conjecture is false for 
simple random walk on a Galton-Watson tree conditioned to survive 
when the offspring distribution is in the domain of attraction of a stable 
law with index $\alpha\in(1,2)$. (Similar result for the annealed case was already  
obtained in \cite{Kest86a}.) 
In both papers, the key ingredient is to obtain suitable bounds on volume and 
effective resistance mentioned above that are sufficient  
for general quenched and annealed estimates. 

The organization of the paper is as follows. 
In the next subsection, we summarize the framework and the main results
on random graphs. In Section 2, we give the application of our 
main results  
to the long range percolation. In Section 3, we give the full proof of the 
main results. 
Although the principal ideas of the proof is quite similar
to the ones in Section 1--2 of \cite{BJKS}, lots of additional careful 
computations are needed to obtain this general version. So, we think it
would help readers to include the full proof.  
\subsection{Framework and Main results}
\label{sec-rwre}

Let $\Gamma = (G,E)$ be an infinite graph, with vertex set $G$
and edge set $E$.  The edges $e \in E$ are \emph{not} oriented.
We assume that $\Gam$ is connected.
We write $x\sim y$ if $\{x,y \} \in E$, and
assume that $(G, E)$ is locally finite,
i.e., $\mu_y<\infty$ for each $y\in G$,  where $\mu_y$ is the number of
bonds that contain $y$. 
Note that $\mu_x\ge 1$ since $\Gamma$ is connected. 
We extend $\mu$ to a measure on $G$.
Let $d(\cdot,\cdot)$ be a metric on $G$. (Note that $d$ is not necessarily a
graph distance. Any metric on $G$ may be used in this section.) We write
\eq
     B(x,r)=\{y: d(x,y) <  r\}, \qquad V(x,r)=\mu(B(x,r)), \quad r \in (0,\infty).
\en
We call $V(x,r)$ the \emph{volume} of the ball $B(x,r)$.
We will assume $G$ contains a marked vertex, which we denote $0$,
and we write
 \eq\label{vrd}
   B(R) = B(0,R), \qquad V(R)= V(0,R).
\en

Let $X =(X_n, n  \in  \bZ_+, P^x, x \in G)$ be the discrete-time
simple random walk on $\Gamma$.  Then $X$ has transition probabilities
\eq
    P^x(X_1=y) = \frac{1}{\mu_x}, \quad y \sim x.
\en
We define the transition  density (or discrete-time heat kernel) of $X$ by
\eq
\lbeq{TAdef11}
     p_n(x,y)= P^x(X_n=y)\frac{1}{\mu_y};
\en
we have $p_n(x,y)=p_n(y,x)$. For $A \subset G$, we write
\eq
\lbeq{TAdef}
     T_A =\inf\{n\ge 0: X_n \in A\}, \qquad \tau_A = T_{A^c},
\en
and let
\eq
    \tau_R = \tau_{B(0,R)} =\min\{ n \ge 0: X_n \not\in B(0,R)\}.
\en
We define a quadratic form $\Ecal$ by
\eq
      \Ecal(f,g)=\tfrac 12\sum_{\substack{x,y\in G \\ x \sim y}}
      (f(x)-f(y))(g(x)-g(y)).
\en
If we regard $\Gamma$ as an electrical network with a unit resistor on each
edge in $E$, then $\Ecal(f,f)$ is the energy dissipation when the
vertices of $G$ are at a potential $f$.
Set $H^2=\{ f\in  \bR^{G}: \Ecal(f,f)<\infty\}$.
Let $A,B$ be disjoint subsets of $G$.
The effective resistance between $A$ and $B$ is defined by:
\begin{equation}
\label{3.3bk}
    \Reff(A,B)^{-1}=\inf\{\Ecal(f,f): f\in H^2, f|_A=1, f|_B=0\}.
\end{equation}
Let $\Reff(x,y)=\Reff(\{x\},\{y\})$, and $\Reff(x,x)=0$.

It is known that $\Reff(\cdot,\cdot)$ is a metric on $G$ (see
\cite[Section~2.3]{Kiga01}), and the following holds.
\begin{equation}\label{ressob}
|f(x)-f(y)|^2\le \Reff(x,y)\Ecal (f,f), \qquad\forall f\in L^2(G,\mu).
\end{equation}

\medskip

We now consider a probability space $(\Omega, \Fcal, \bP)$ carrying a
family of random graphs $\Gam(\omega)=(G(\omega),E(\omega), \omega\in \Omega)$.
We assume
that, for each $\omega \in \Omega$, the graph $\Gam(\omega)$ is infinite,
locally finite, connected, and contains a marked vertex $0\in G$.
Let $d(\cdot,\cdot):=d_\omega(\cdot,\cdot)$ be a metric on $G(\omega)$. 
(Again, $d$ is any metric on $G(\omega)$, not necessarily the graph distance.)   
We denote balls in $\Gam(\omega)$ by  $B_\omega(x,r)$, their volume
by $V_\omega(x,r)$, and write
\eq
    B(R) = B_\omega(R) = B_\omega(0,R),
    \quad\quad
    V(R) = V_\omega(R) = V_\omega(0,R).
\en
We write
$X=(X_n , n \ge 0, P_\omega^x, x \in G(\omega))$ for the
simple random walk on $\Gam(\omega)$, and denote
by $p_n^\omega(x,y)$ its transition density with respect
to $\mu(\omega)$.
To define $X$ we introduce a second
measure space $(\overline \Omega, \overline \Fcal)$,
and define $X$ on the product $\Omega \times \overline \Omega$.
We write
$\overline \omega$ to denote elements of $\overline \Omega$.

Let $v, r: \bN\to [0,\infty)$ be strictly increasing
functions with $v(1)=r(1)=1$ which satisfy 
\begin{equation}\label{vrdoubl}
C_1^{-1}\Big(\frac R{R'}\Big)^{d_1}\le \frac{v(R)}{v(R')}\le C_1
\Big(\frac R{R'}\Big)^{d_2},~~
C_2^{-1}\Big(\frac R{R'}\Big)^{\al_1}\le \frac{r(R)}{r(R')}\le C_2
\Big(\frac R{R'}\Big)^{\al_2}
\end{equation}
for all $0<R'\le R<\infty$, where $C_1,C_2\ge 1$, $1\le d_1\le d_2$ 
and $0<\al_1\le \al_2\le 1$. For convenience, we let $v(0)=r(0)=0$,
$v(\infty)=r(\infty)=\infty$ and extend them to 
$v, r: [0,\infty]\to [0,\infty]$  
such that $v,r$ are continuous, 
strictly increasing, and satisfy (\ref{vrdoubl}). One such extension is
to extend them linearly. 

The key ingredients in our analysis of the simple random walk are volume
and resistance bounds.  The following defines a random set $J(\lambda)$ 
of values of $R$ for which we have `good' volume and effective resistance
estimates.

\begin{defn}\label{jdef}
\emph{ Let $\Gamma=(G,E)$ be as above.
For $\lambda>1$, define
\begin{eqnarray*}
J(\lambda)&=&\{R \in [1,\infty]: \lambda^{-1} v(R)\le V(R) \le \lambda v(R),
\Reff (0, B(R)^c)\ge \lambda^{-1} r(R),\\
&&~~ \Reff (0,y)\le \lam r(d(0,y)),~\forall
y\in B(R)\}.
\end{eqnarray*}
}\end{defn}

As we see, $v(\cdot)$ gives the volume growth order and $r(\cdot)$  
gives the resistance growth order.

We now make the following assumptions concerning the graphs $(\Gam(\omega))$.
This involves upper and lower bounds on the volume, as well as an estimate
which says that $R$ is likely to be in $J(\lambda)$ for large enough $\lambda$. 
Note that some assumption includes another assumption, because we assume
part of them for each theorem and proposition.

\begin{ass}\label{ass-rwre}
(1) There exist $\lambda_0>1$ and $p(\lambda)$ which goes to $0$ as 
$\lambda\to \infty$ such that,
\begin{equation}\label{plamassump}
\bP (R\in J(\lambda))\ge 1-p(\lambda)\qquad \mbox{for each}~~R\ge 1, \lambda\ge\lam_0.
\end{equation}
(2) $\bE [\Reff (0,B(R)^c)V(R)]\le c_{1} v(R)r(R).$\\
(3) There exist $q_0, c_2>0$ such that
\eq
  p(\lam) \le \frac{c_2}{\lam^{q_0}}.
\en
(4) There exist $q_0, c_3>0$ such that
\eq
  p(\lam) \le \exp (-c_3\lam^{q_0}).
\en
\end{ass}

We have the following consequences of
Assumption~\ref{ass-rwre} for random graphs. 
Some of the results apply also to
the random walk started from an arbitrary point $x \in G(\omega)$.
Some statements in the first proposition involve the annealed law
\eq
    P^* = \bP \times P^0_\omega.
\en

Let ${\cal I}(\cdot)$ be the inverse function of $(v\cdot r)(\cdot)$. 
\begin{prop} \label{ptight}
Suppose Assumption \ref{ass-rwre}(1) holds.
 Let $n\ge 1$, $R\ge 1$. Then
\begin{align}
\label{pt-a}
 \bP (\theta^{-1}\le \frac {E^0_\omega\tau_R}{v(R)r(R)}\le \theta) &\to 1
\quad \text{ as } \theta\to \infty, \\
\label{pt-b}
  \bP (\theta^{-1}\le v({\cal I}(n))p_{2n}^\omega(0,0)\le \theta) &\to 1
\quad \text{ as } \theta\to \infty, \\
\label{pt-c}
  P^*( \frac{d(0,X_n)}{{\cal I}(n)} < \theta) &\to 1
\quad \text{ as } \theta\to \infty. \\
\label{pt-d}
 P^*( \theta^{-1} < \frac{1+d(0,X_n)}{{\cal I}(n)} ) &\to 1
\quad \text{ as } \theta\to \infty.
\end{align}
In each case the convergence is 
uniform in $n$ and $R$.
\end{prop}

\smallskip
Since $P^0_\omega(X_{2n}=0) \approx 1/v({\cal I}(n))$, we cannot replace
$1+d(0,X_n)$ by $d(0,X_n)$ in (\ref{pt-d}).

\begin{prop} \label{pmeans}
Suppose Assumption 
\ref{ass-rwre}(1) and (2) hold. Then
\eqalign
\label{e:mmean}
    c_1v(R)r(R) \le \bE (E^0_\omega\tau_R) \le c_2v(R)r(R) & \text{  for all } R\ge 1,
\\
\label{e:pmean}
 \frac{c_3}{v({\cal I}(n))}\le \bE (p_{2n}^\omega(0,0))
 &\text{  for all } n \ge 1,
\\
\label{e:dmean}
       c_4 {\cal I}(n)\le \bE (E_\omega^0 d(0,X_n)) &\text{  for all } n\ge 1.
\enalign
Assume in addition that 
there exist $c_5>0, \lam_0>1$ and $q_0'>2$ such that 
\begin{equation}\label{newcond45}
\bP (\lambda^{-1} v(R)\le V(R),~ \Reff (0,y)\le 
\lam r(d(0,y)),~\forall y\in B(R))\ge 1-\frac {c_5}{\lam^{q_0'}},
\end{equation}
for each $R\ge 1, \lambda\ge\lam_0$. 
Then 
\begin{equation}
\label{e:pmeanII}
\bE (p_{2n}^\omega(0,0))\le 
 \frac{c_6}{v({\cal I}(n))}\text{  for all } n \ge 1.
 \end{equation}

\end{prop}

We do not have an upper bound in (\ref{e:dmean}); see 
Example~2.6 in \cite{BJKS}. 

The additional assumption that $p(\lam)$ decays 
either polynomially or exponentially enable us to obtain limit theorems.
Both of the following theorems refer to the random walk started at an
arbitrary point $x \in G(\omega)$.

\begin{theorem}
\label{thm-rwre} (I) Suppose Assumption \ref{ass-rwre}(1) and (3) hold.
Then there exist $\beta_1, \beta_2, \beta_3, \beta_4 < \infty$, and a subset $\Omega_0$
with $\bP(\Omega_0)=1$ such that the following statements hold.\\
(a) For each $\omega \in \Omega_0$ and $x \in G(\omega)$ there
exists $N_x(\omega)< \infty$ such that
\eq
 \label{e:logpnlima}
   \frac{(\log n)^{-\beta_1}}{v({\cal I}(n))} \le
 p^\omega_{2n}(x,x) \le \frac{(\log n)^{\beta_1}}{v({\cal I}(n))}, \quad  n\ge N_x(\omega).
\en
(b)  For each $\omega \in \Omega_0$ and $x \in G(\omega)$ there
exists $R_x(\omega)< \infty$ such that
\eq
 \label{e:logtaulima}
   (\log R)^{-\beta_2} v(R)r(R) \le E^x_\omega \tau_R
 \le  (\log R)^{\beta_2} v(R)r(R), \quad  R\ge R_x(\omega).
\en
(c)  Let $Y_n= \max_{0\le k \le n} d(0,X_k)$.
For each $\omega \in \Omega_0$ and $x \in G(\omega)$
there exist $N_x(\omega,\overline \omega), R_x(\omega, \overline \omega)$
such that $P^x_\omega( N_x <\infty)=P^x_\omega( R_x <\infty)=1$,
and such that
\begin{align}
 \label{e:ynlim}
 (\log n)^{-\beta_3} {\cal I}(n) \le Y_n(\omega, \overline \omega)
 &\le (\log n)^{\beta_3} {\cal I}(n),
 \quad n \geq N_x(\omega, \overline \omega), \\
   \label{e:rnlim}
 (\log R)^{-\beta_4} v(R)r(R) \le \tau_R(\omega, \overline \omega)
 &\le (\log R)^{\beta_4}  v(R)r(R), \quad\quad
 R \geq R_x(\omega, \overline \omega).
\end{align}
(II) Suppose Assumption \ref{ass-rwre}(1) and (4) hold.
Then there exist $\beta_1, \beta_2 < \infty$, and $\Omega_0$
with $\bP(\Omega_0)=1$ such that 
(\ref{e:logpnlima}) and (\ref{e:logtaulima}) hold with $\log\log n$ (resp. 
$\log\log R$) instead of $\log n$ (resp. $\log R$).\\
(III)  Suppose Assumption \ref{ass-rwre}(1) and (3) hold. Suppose further 
that $v, r$ satisfy the following in addition to (\ref{vrdoubl}): 
\eq\label{logvrcond}
C_3^{-1}R^D(\log R)^{-m_1}\le v(R)\le C_3R^D(\log R)^{m_1},~~
C_4^{-1}R^\al(\log R)^{-m_2}\le r(R)\le C_4R^\al(\log R)^{m_2}
\en
for all $0<R$, where $C_3,C_4\ge 1$, $1\le D$, $0<\al\le 1$
and $0<m_1,m_2$. Then the following statements hold.\\
(a) $d_s(G) :=-2\lim_{n\to\infty}
\frac{\log p_{2n}^\omega(x,x)}{\log n}= \frac{2D}{D+\al}$, $\bP$--a.s.,
and the random walk is recurrent.\\
(b) $\lim_{R \to \infty} \frac{\log  E^x_\omega \tau_R}{\log R} =D+\al$.\\
(c) Let $W_n = \{X_0,X_1,\ldots, X_n\}$ and let
$S_n = \mu(W_n) = \sum_{x \in W_n} \mu_x$.
For each $\omega \in \Omega_0$ and $x \in G(\omega)$,
\eq
\label{e:snlim}
 \lim_{n \to \infty} \frac{\log S_n}{\log n} = \frac D{D+\alpha}, \quad
 P^x_\omega \text{--a.s.}.
\en
\end{theorem}
\begin{rem}\label{resremvi}
1.  Let
\begin{eqnarray*}
{\hat J}(\lambda)&=&\{R \in [1,\infty]: \lambda^{-1} v(R)\le V(R) \le \lambda v(R),
\Reff (0, B(R)^c)\ge \lambda^{-1} r(R),\\
&&~~ \Reff (0,y)\le  r(d(0,y)),~\forall y\in B(R)\}.
\end{eqnarray*}
(Note that $J(\lambda)$ contains $\lambda$ in the last inequality
whereas ${\hat J}(\lambda)$ does not.) Assume that Assumption \ref{ass-rwre}(1)
holds w.r.t. ${\hat J}(\lambda)$ and further the following holds:\\
\eq\label{e:addinvvol}
\bE [1/V (R)]\le c_{1}/v(R).\en 
(Note that this condition is a bit weaker than 
(\ref{newcond45}).)    
Then, (\ref{e:pmean}) and (\ref{e:pmeanII}) hold. 

\noindent 
2. If one chooses the resistance metric $\Reff (\cdot,\cdot)$ as the 
metric $d(\cdot,\cdot)$, then clearly $\Reff (0,y)\le \lam r(d(0,y))$ 
holds with $r(x)=x$ and $\lam=1$. \newline
3. If all the vertices in $G(\omega)$ have degree bounded by a constant
$c_0$ for all $\omega\in \Omega_0$, then $|W_n| \le S_n \le c_0 |W_n|$. Hence,
under the assumption (\ref{logvrcond}),  
Theorem \ref{thm-rwre}(III)(c) implies also that
\eq
\label{e:oaocpnj}
\lim_{n \to \infty} \frac{\log |W_n|}{\log n} = \frac D{D+\alpha}, \quad
 P^x_\omega \text{--a.s.}
 \en
\end{rem}

\section{Application: Long range percolation}
In this section, we will apply the theorem in Section 1
to the long range percolation.
Let $\mathbf{p}={\{p(n)\}}_{n=1}^{\infty}$ be a sequence of real numbers.
Each unoriented pair of distinct points $x,y\in {\mathbb{Z}}^d$ 
is connected by an unoriented bond with 
probability $p(x,y)=p(y,x)=p(|x-y|)$, independently of other pairs.
Here, $|x-y|=\sum_{i=1}^{d} |{x_i}-{y_i}|$.
We consider the situation that $\mathbf{p}$ satisfies
\eq\label{con-prob}
\lim_{n\to\infty}\frac{p(n)}{\beta n^{-s}}=1,
\en  
for some $s>0$, $\beta >0$.

Let ${\mu}_{xy}$ be $\{0,1\}$-valued random variable, 
which takes $1$ if $x$ and $y$ are connected by a bond 
and takes $0$ if there is no bond between $x,y$.
(${\mu}_{xy}={\mu}_{yx}$, and ${\mu}_{xx}=0$.) 
${\mu}_x =\sum_{y\in {\mathbb{Z}}^d}{\mu}_{xy}$ stands for the number 
of bonds which have $x$ as an endpoint.
$G={\mathbb{Z}}^d$ is the vertex set and $E=\{\langle
x,y\rangle|{\mu}_{xy}=1\}$
is the edge set of the corresponding random graph. 
(We identify $\langle x,y\rangle =\langle y,x\rangle$.)

\bigskip

Here, we give some comments on the backgrounds.
Random walks on long-range percolation clusters 
is discussed in \cite{Ber02}.
Let $\mathbf{p}$ be a sequence satisfying (\ref{con-prob}) and 
$p(n)\in[0,1)$ for $n\ge 1$.
The random graphs are 
locally finite if and only if $s>d$, and we can define random walks in such a case. 
We choose $\mathbf{p}$ for which  
there exists a unique $\infty$-cluster with probability $1$.
Then, the main results in \cite{Ber02} are the following: 

\begin{flushleft}
(1) For $d=1$, random walks are 
transient if $1<s<2$, and recurrent if $s=2$.

(2) For $d=2$, random walks are 
transient if $2<s<4$, and recurrent if $s\ge 4$.
\end{flushleft}

In the above, 
the case $d=1$, $s>2$ is not mentioned because 
there is no $\infty$-cluster in such a case.
From now on, we explore the case $d=1$, 
$p(1)=1$, and for $n\ge 2$, $p(n)\in [0,1)$ satisfies (\ref{con-prob})
for some $s>2$. In this case, the effects of long bonds are not 
so strong, and as we will see later, behaviours of random walks 
are similar to the random walk on $\bZ$. Also, we will refer to some
kind of discontinuity on $s=2$.  

Let $d(x,y):=|x-y|$ and define $B(R), V(R)$ with respect to this
metric. Then, since $V(R)\ge \sharp B(R)=2R-1$, 
\[\bP (V(R)\le \lam^{-1}R)=0\qquad\text{if }~~\lam>1.\]
Further, we have
\[\bP (V(R)\ge \lam R)\le \frac{\bE V(R)}{\lam R}\le 
\frac c{\lam}.\]
The upper bound on the resistance is obvious by comparing 
percolation clusters with $\bZ$ ;
\eq\label{res-ub}
\bP\left[\bigcup_{y\in B(R)}\{R_{\mathrm{eff}}(0,y)>\lam d(0,y)
\}\right]=0\qquad\text{if }~~\lam>1,
\en
and
\eq\label{res-meanub}
\bE [\Reff (0,B(R)^c)V(R)]\le R\bE[V(R)]\le c_{1} R^{2}.
\en
The remained work 
is the lower bound on the resistance. We have the following.
\begin{prop}\label{res-1}
Let $q=1$ for $s>3$, and let $q$ be an any value in $(0, s-2)$ for 
$2<s\le 3$.
Then, there exists $c_1=c_1(\beta, s, q)>0$, such that for each $R\ge 1$,
\begin{equation}\label{res-lb}
\bP[\Reff (0, {B(R)}^c)<{\lam}^{-1}R]
\le c_1 {\lambda}^{-q}.
\end{equation}
\end{prop}
\proof
First, we apply the ``projecting long bonds'' method in Lemma 
3.8 in  
\cite{Ber02} to our case.
For each $\omega$, we construct a new weighted graph from the original
one in the following way.

\begin{flushleft}
(1) If a bond $\langle x,y\rangle$ 
such that $x,y\in \bZ$, $x+2\le y$ exists, then,
divide the bond into $y-x$ short bonds.

(2) For each $i=1,\cdots ,y-x$, replace the $i$-th short bond by 
a bond which has $x+i-1$ and $x+i$ as its endpoints and has 
weight $y-x$. 

(3) Repeat (1),(2) for all bonds except nearest-neighbour bonds.  
\end{flushleft}
We use the notation $\tReff$ for the resistance on the new graph. 
By the \emph{shorting law}  
in the terminology of the electrical 
network, the resistance does not increase in the above procedures.
And from the way of construction, we can see that   
\begin{equation}
\tReff (0,R)= \sum_{i=1}^{R}(\sum_{e\in A_{i}}|e|)^{-1},
\end{equation}
where $A_i$ is the set of all $\langle u,v
 \rangle$ satisfying $u,v\in \bZ$, $u<v$, ${\mu}_{uv}=1$, and
 $[u,v]\supset[i-1,i]$.
(In other words, $A_i$ is the collection of bonds crossing over
 $[i-1,i]$.) 
We denote $|e|=|u-v|$ for $e=\langle u,v \rangle$. 
It is easy to see that
\begin{eqnarray*}
\bE[{\Reff (0, {B(R)}^c)}^{-q}]&\le&\bE[{\tReff(0,
 {B(R)}^c)}^{-q}]\\
&=&\bE[\{{\tReff(0,R)}^{-1}+
{\tReff(0,-R)}^{-1}
\}^q]\\
&\le&\bE[{\tReff(0,R)}^{-q}+
{\tReff(0,-R)}^{-q}]=2\bE[{\tReff(0,R)}^{-q}],
\end{eqnarray*}
and
\begin{eqnarray*}
\bE[{\tReff (0,R)}^{-q}]&=&
\bE\left[\{\sum_{i=1}^{R}(\sum_{e\in A_{i}}|e|)^{-1}
\}^{-q}\right]\\
&\le& R^{-q-1}\sum_{i=1}^{R}\bE\left[(\sum_{e\in A_{i}}|e|
)^{q}\right].
\end{eqnarray*}
We have used the H\"{o}lder inequality in the last estimate.
The expectation in the right hand side is finite 
for each $q$ ; 
\begin{eqnarray*}
\bE\left[(\sum_{e\in A_{i}}|e|)^{q}\right]&\le&
\bE[\sum_{e\in A_{i}} {|e|}^{q}]\\
&=&\bE[\sum_{n=1}^{\infty}\sum_{k=1}^{n}1_{\{{\mu}_{k-n,k}=1\}} n^q]
\\
&=& \sum_{n=1}^{\infty} n^q 
\sum_{k=1}^{n}\bP[{\mu}_{k-n,k} =1]\\
&\le&c_2 \sum_{n=1}^{\infty}n^{q-s+1} < \infty .
\end{eqnarray*}
Combining these calculations, we have
\begin{equation}\label{res-q}
\bE[{\Reff(0, {B(R)}^c)}^{-q}]\le c_3 R^{-q}.
\end{equation}
So, by the Chebyshev inequality and (\ref{res-q}), 
\begin{eqnarray*}
\bP[\Reff (0, {B(R)}^c)\le {\lam}^{-1}R]&\le&
{\lam}^{-q} R^q \bE[{\Reff(0, {B(R)}^c)}^{-q}]\\
&\le&c_3 {\lam}^{-q},
\end{eqnarray*}
which completes the proof.
\qed

\bigskip
From the above estimates, we have the following. Below,
$a_n\sim b_n$ stands for $\lim_{n\to\infty}a_n/b_n=1$. 
\begin{theorem}\label{lr-result}
The long-range percolation on $\bZ$ with 
$p(n)\sim \beta n^{-s}$ for 
$s>2$, $\beta >0$, $p(1)=1$
satisfies Assumption 
\ref{ass-rwre}(1),(2) and (3) 
with $v(x)=r(x)=x$.
\end{theorem}
In this case, the additional assumption (\ref{newcond45}) 
also holds directly from the condition $p(1)=1$. So, 
we obtain the conclusion of Proposition \ref{ptight}, \ref{pmeans} 
and Theorem \ref{thm-rwre}(I,III) with $v(x)=r(x)=x$.

\begin{rem}\label{lr-rem}
(1) 
We have proved that, when $s>2$, 
$p_{2n}(x,x)$ is the order of $n^{-\frac{1}{2}}$. 
On the other hand, from the transience result in \cite{Ber02}, 
it is natural to see that, when $1<s<2$,  
$p_{2n}(x,x) \sim n^{-\xi(s)}$ 
($\xi(s) >1$) in some sense, though we do not have a rigorous proof.   
Hence, Theorem \ref{lr-result} implies that
the order of the heat kernel is \emph{discontinuous} at $s=2$.

(2)
In the one-dimensional long-range percolation model, 
the phenomena at the point $s=2$ are non-trivial. 
In \cite{AN86},  
the discontinuity of the percolation density 
at $s=2$ is shown. 
In the recent study in \cite{BBY07} , 
long range percolation mixing time is considered and 
it is shown that
the order of the mixing time changes discontinuously when $s=2$.  
In \cite{Mis07}, estimates of effective resistance are given 
for general $d$, $s$, and the discontinuity mentioned in (1) 
is shown in a sense of effective resistance. 

(3)
Note that the proof of Theorem \ref{lr-result} can be also applied
to some case where we do not assume the independence of open bonds. 
$s=2$ is not necessarily critical in the case.
\end{rem}

\begin{rem}\label{lr-rem2}
There is a question whether the log-corrections in Theorem
\ref{thm-rwre}(I) can be weaken or not. To answer this question, it is 
crucial to study the fluctuation of the volume and resistance.  
First, let us consider the fluctuation of the volume.
The following large deviation estimate holds for $s>1$.
\[
\bP[V(R)\ge c_1 R]\le \exp\{-c_2 R\}.
\]
Since $\sum_{R=1}^{\infty} \exp\{-c_2 R\} <\infty$,  
by the Borel--Cantelli lemma, we have 
\[
\bP[\limsup _{R\to\infty} \frac{V(R)}{R}\le c_1]=1. 
\]
The lower bound of $V(R)/R$ is trivial from $p(1)=1$.
Therefore, there is no fluctuation of the volume such as  
Proposition $2.8$ in \cite{BK06}.

Next, we consider the fluctuation of the resistance.
When $2<s<3$, by calculating $\bE [\Ecal (g,g)]$ for appropriate $g$,
we can see that $\bE [{\Reff (0, {B(R)}^c )}^{-1}]
\le c R^{2-s}$. It seems to us that this is the best estimate one can
obtain. If so, there may be some fluctuation of the resistance. 
\end{rem}
We also give an example satisfying Assumption \ref{ass-rwre} (4) 
instead of (3).

\begin{theorem}\label{exp-case}
We consider 
the long-range percolation on $\bZ$ with $p(n)\sim e^{-cn}$
for $c>0$, $p(1)=1$. Then, 
Assumption 
\ref{ass-rwre}(1),(2) and (4) are 
satisfied with $v(x)=r(x)=x$. 
\end{theorem}

\proof
The estimate for volume is easy, and the upper bound 
of resistance is trivial.
We will show that $\bP[\Reff (0, {B(R)}^c)\le {\lam}^{-1}R]\le 
e^{-{c_1}\lam}$ for some ${c_1}>0$.
\begin{eqnarray*}
\bP[\Reff (0, {B(R)}^c)\le {\lam}^{-1}R]
&=&\bP[\inf \{\Ecal (f,f) : f(0)=1, f|_{{B(R)}^c} \equiv 0\}
\ge \lam R^{-1}]\\
&\le&\bP[\Ecal (g,g)\ge \lam R^{-1}]\\
&\le& e^{-{c_2}\lam} \bE[\exp\{{c_2}R \Ecal (g,g)\}],
\end{eqnarray*}
where $g(x)=1_{B(R)}(x) (1-{R^{-1}}|x|)$, and 
$c_2$ is a positive constant determined later. We see that
\begin{eqnarray*}
\bE[\exp\{{c_2}R \Ecal (g,g)\}]&=&
\bE[\exp\{{c_3}R \sum_{x,y\in\bZ}{|g(x)-g(y)|}^2 {\mu}_{xy}\}]\\
&=& \prod_{x,y\in\bZ} \bE [\exp\{{c_3}R {|g(x)-g(y)|}^2 
{\mu}_{xy}\}]\\
&=& \prod_{x,y\in\bZ} \left\{1+(\exp\{c_3 R {|g(x)-g(y)|}^2
\}-1) e^{-c|x-y|}
\right\}
\equiv \prod_{x,y\in\bZ} I_{xy}, 
\end{eqnarray*}
where $c_3 =\frac{c_2}{2}$.
Clearly, $\prod_{x,y\in {B(R)}^c} I_{xy} =1$, and
\begin{eqnarray*}
\prod_{x,y\in B(R)} I_{xy} &\le&
\prod_{x,y\in B(R)}\left\{1+(\exp\{c_3 R^{-1} {|x-y|}^2
\}-1) e^{-c|x-y|}
\right\}\\
&\le& \prod_{n=1}^{2R} \{1+(\exp(c_3 R^{-1} n^2 ) -1) e^{-cn}
\}^{2R} \equiv X. 
\end{eqnarray*}
Now 
\begin{eqnarray*}
\log X &=& \sum_{n=1}^{2R} 2R\log \{1+(\exp(c_3 R^{-1} n^2 ) 
-1) e^{-cn} \}\\
&\le& 2R \sum_{n=1}^{2R} (\exp(c_3 R^{-1} n^2 ) 
-1) e^{-cn}\\
&=& 2R \{\sum_{n=1}^{[\sqrt{R}]} +\sum_{n=[\sqrt{R}] +1}^{2R} \}
\equiv 2R(S_1 + S_2).
\end{eqnarray*}
We first estimate $S_1$. 
\[S_1 \le \sum_{n=1}^{[\sqrt{R}]} c_4 R^{-1} n^2 e^{-cn}
\le c_4 R^{-1} \sum_{n=1}^{\infty} n^2 e^{-cn} \le c_5 R^{-1}.\]
For the estimate of $S_2$, we choose $c_3$ sufficiently small so that  
\[S_2 \le \sum_{n=[\sqrt{R}] +1}^{2R} \exp\{c_3 R^{-1} n^2
\} e^{-cn}
\le \sum_{n=[\sqrt{R}] +1}^{2R} \exp\{-(c-2{c_3})n\}
\le e^{-c_6 \sqrt{R}}. \]
Therefore, we have $\log X \le c_7$ , and $X\le e^{c_7}$.
Furthermore, we see that
\begin{eqnarray*}
\prod_{x\in B(R), y\in {B(R)}^c} I_{xy} &=&
\prod_{x\in B(R)} \prod_{y\in {B(R)}^c} 
\left\{1+(\exp\{c_3 R^{-1} {(R-|x|)}^2
\} -1)e^{-c|x-y|}\right\}\\
&\equiv& \prod_{x\in B(R)} a_x ,
\end{eqnarray*}
and
\begin{eqnarray*}
\log a_x &=& \sum_{y\in {B(R)}^c } \log \left\{1+(\exp\{c_3 R^{-1} {(R-|x|)}^2
\} -1)e^{-c|x-y|}\right\}\\
&\le& \sum_{y\in {B(R)}^c } (\exp\{c_3 R^{-1} {(R-|x|)}^2
\} -1)e^{-c|x-y|}\\
&\le& c_8 \exp\{c_3 R^{-1} {(R-|x|)}^2
\} e^{-c(R-|x|)}\\
&\le& c_8 e^{c_3 (R-|x|)} e^{-c(R-|x|)}
= c_8 e^{- c_9 (R-|x|) },
\end{eqnarray*}
for sufficiently small $c_3$. Thus,
\begin{eqnarray*}
\prod_{x\in B(R)} a_x &\le& \prod_{x\in B(R)}
\exp\{c_8 e^{- c_9 (R-|x|) } \} \\
&=& \exp\{c_8 \sum_{x\in B(R)} e^{- c_9 (R-|x|) } \}
\le \exp\{2 c_8 \sum_{n=1}^{\infty} e^{- c_9 n} \}
\le c_{10}.
\end{eqnarray*}
From these estimates, we obtain the result for suitable $c_1$ .
\qed

\section{Proof of the main results}
\label{sec:rwresults}

In Section~\ref{sec-gg}, we prove several preliminary results for random
walk on a fixed but general graph.  Then in Section~\ref{sec-genrwre} we
apply these results to prove Propositions~\ref{ptight}--\ref{pmeans} and
Theorem~\ref{thm-rwre}. 
We adopt the convention that if we cite elsewhere the constant $c_1$
in Lemma $3.2$ (for example), we denote it as $c_{3.2.1}$.
$C_1,C_2$ stand for the constants in (\ref{vrdoubl}). 

\subsection{Estimates for general graphs}
\label{sec-gg}

In this section, we fix an infinite locally-finite connected
graph $\Gam = (G,E)$, and
use bounds on the quantities $V(R)$ and $\Reff(0, B(R)^c)$ to
control $E^0 \tau_R$, $p_n(0,0)$ and $E^0 d(0,X_n)$, where 
$d(\cdot,\cdot)$ is a metric on $G$. 
To deal with issues related to the possible
bipartite structure of the graph, we will consider
$p_n(x,y)+p_{n+1}(x,y)$.

\begin{prop}\label{nash}
Let $0\in G$ and $f_n(y)=p_n(0,y)+p_{n+1}(0,y)$. \\
(a) Let $R\ge 1$ and assume that 
\begin{equation}\label{vrcond}
V(0,R)\ge \lam^{-1}v(R),~\Reff (0,y)\le \lam r(d(0,y))\qquad
\mbox{for all }~ y\in B(R). 
\end{equation}
Then
$$ f_{n}(0) \le \frac{c_1 \lam}{v({\cal I}(n))}\qquad
\mbox{for}~~\frac 12 v(R)r(R)\le n\le v(R)r(R). $$
(b) We have
$$ | f_n(y)-f_n(0) |^2
 \le  \frac {c_2}{n} \Reff (0,y) p_{ 2 \lfloor n/2 \rfloor} (0,0) . $$
\end{prop}
\proof (a) A natural modification of 
the third equation in \cite[Proposition~3.3]{BCK05} using 
$\Reff (0,y)\le  \lam r(d(0,y))\le \lam r(R)$ gives
\[ f_{2n}(0)^2 \le \frac{c}{V(0,R)^2} + \frac{c \lam r(R) f_{2n}(0)}{n}, 
\qquad  \mbox{ for all }  n\ge 1,~R>0. \]
Using $a+b \leq 2(a \vee b)$, we see that
$f_{2n}(0) \le  ({c'}/{V(0,R)}) \vee ({c' \lam r(R)}/{n})$. So, by setting 
$v(R)r(R)/2\le n\le v(R)r(R)$, 
\[
f_{2n}(0) \le (c'/V(0,R)) \vee (2c'\lam/v(R))\le 2c'\lam/v(R)
\le c_1\lam/v({\cal I}(n)), \]
where we used $V(0,R)\ge \lam^{-1}v(R)$ in the second inequality. 

(b) Using (\ref{ressob}), 
\[ |f_n(y)-f_n(0)|^2 \le \Reff(0,y) \Ecal(f_n,f_n). \]
We then use \cite[Lemma~3.10]{BCK05} to bound $\Ecal(f_n,f_n)$. \qed

\begin{prop}\label{bkb}
Let $R\ge 1$, $m \ge 1$, 
$\eps^{\al_1} \le 1/(4C_2m\lam)$. 
Write $B=B(0,R)$, $B'=B(0, \frac12\eps R)$,
$V=V(0,R)$, $V'=V(0,\frac12\eps R)$ and suppose $\Reff(0,y)\le \lam r(d(0,y))$
for all $y\in B(R)$.\\
(a) For $x\in B$,
\eq\label{e:dsad}
    E^x \tau_R \le 2\lam r(R) V.
\en
(b) Suppose further that
\eq
\label{e:relb}
\Reff(x, B^c) \ge r(R)/m  \quad \text{ for } x \in B(0,\eps R).
\en
Then for $x \in  B'$,
\begin{align}
\label{e:tmlb}
   E^x \tau_R &\ge  \frac{r(R) V'}{4m}, \\
\label{e:tplb}
   P^x( \tau_R > n) &\ge \frac{V'}{8 m\lam V} - \frac{n}{2 r(R) \lam V} \quad
\text{ for } n \ge 0, \\
\label{e:plb}
 p_{2n}(x,x) &\ge \frac{c_1 (V')^2}{m^2 \lam^2 V^{3} } \quad \quad
\text{ for } n \le  \frac{r(R) V'}{8 m }.
\end{align}
\end{prop}

\proof
For any $z\in B$ we have 
$$  E^z\tau_B = \sum_{y\in B} g_B(z,y)\mu_y, $$
where $g_B(\cdot,\cdot)$ is the Green kernel of the Markov chain killed on exiting 
$B$. 
 
(a) Since $\Reff(z,B^c)\le \Reff (0,z)+\Reff (0,B^c)\le 2\lam r(R)$ 
for any $z\in B$,  
\eq\label{e:newads}
   E^z\tau_B= \sum_{y\in B} g_B(z,y)\mu_y 
 \le \sum_{y\in B} g_B(z,z) \mu_y
  = \Reff(z,B^c) V(0,R) \le 2 \lam r(R) V(0,R), \en
where we used the fact $g_B(z,z)=\Reff(z,B^c)$ in the second equality. 
(For the proof of this fact, see, for example, section $3.2$ in \cite{BCK05}.)

 (b) Let $p^x_B(y)=g_B(x,y)/g_B(x,x)$. Since 
$\Ecal (p^x_B,p^x_B)= \Reff (x,B^c)^{-1}=g_B(x,x)^{-1}$ and so if  $x,y\in B'$
$$ | 1-p^x_B(y)|^2 \le \Reff(x,y) \Reff(x,B^c)^{-1} \le \frac{2\lam r(\eps R)m}{r(R)}
\le 2C_2m \eps^{\al_1}\lam \le 1/2, $$
where (\ref{vrdoubl}) is used in the third inequality. 
Hence $p^x_B(y)\ge 1 - 2^{-1/2}\ge \frac 14$.
So, 
$$  E^x \tau_R \ge \sum_{y\in B'} g_B(x,x)p^x_B(y)\mu_y 
  \ge \frac 14 \mu(B') \Reff(x,B^c) 
  \ge r(R) \mu(B')/(4m). $$
By the Markov property, (\ref{e:dsad}) and (\ref{e:tmlb}), for $x\in B'$,
$$ \frac{r(R) V'}{4m} \le 
E^{x}[\tau_R] \le n+ E^x[1_{\{\tau_R> n\} } E^{X_n}(\tau_R)]
\le n + 2 \lam r(R) V P^x( \tau_R>n),$$
for all $n\ge 1$. Rearranging this gives (\ref{e:tplb}).

By (\ref{e:tplb}), 
$$P^x( X_n \in B) \ge  P^x(\tau_R>n) \ge \frac{ (r(R) V'/4m)-n}{2\lam r(R) V}. $$
So, if $n \le r(R) V'/(8m)$ then
\eq\label{e:mkwl}  P^x( X_n \in B) \ge \frac{c_2 V'}{m\lam V}.  
\en
By the Chapman--Kolmogorov equation and the Cauchy-Schwarz inequality, 
$$  P^x( X_n \in B)^2 = ( \sum_{y \in B} p_n(x,y) \mu_y)^2
\le \mu(B) \sum_{y \in B} p_n(x,y)^2 \mu_y  \le p_{2n}(x,x) V, $$
and using (\ref{e:mkwl}) gives (\ref{e:plb}). \qed

\smallskip

Recall the set $J(\lam)$ defined in Definition~\ref{jdef}.
We will need to know that bounds in the following are polynomial in $\lambda$.
To indicate this, we write $c_i(\lam)$ to denote
positive constants of the form $c_i(\lam) = c_i \lam^{\pm q_i}$.
The sign of $q_i$ is such that statements become weaker as $\lam$
increases.
The following proposition controls the mean escape times and transition
probabilities.

\begin{prop}\label{eqest} Let $\lambda>1$.\\
(1) Suppose that $R \in J(\lam)$. Then there exists $c_1(\lam)$
such that
\begin{align}
\label{e:24-11}
E^z \tau_R &\le 2 \lambda^2 v(R)r(R)  \qquad \quad\text{ for } z\in B(R),\\
\label{e:24-12}
p_{n}(0,0) +  p_{n+1}(0,0) &\le \frac{c_1(\lam)}{v({\cal I}(n))} \quad \text{ if }~
\frac 12 v(R)r(R)\le n\le v(R)r(R),\\
\label{e:24-13}
p_{n}(0,y)+p_{n+1}(0,y)  &\le \frac{c_1(\lam)}{v({\cal I}(n))}  \quad\text{ for }
 y \in B(R) \text{ if }~\frac 12 v(R)r(R)\le n\le v(R)r(R).
\end{align}
(2)  There exist $c_2(\lam),\cdots, c_7(\lam)$ such that, if 
 $R, c_2(\lam)R \in J(\lam)$ , then 
\begin{align}
\label{e:24-21}
 c_3(\lam) v(R)r(R) &\le E^x \tau_R \quad \text{ for } x \in 
 B(c_2(\lam)R), \\
\label{e:24-22}
P^0(\tau_R>c_4(\lam)v(R)r(R))&\ge c_5(\lam),\\
\label{e:24-23}
p_{2n}(0,0) &\ge \frac{c_6(\lam)}{v({\cal I}(n))} \quad
\text{ for } \tfrac12 c_7(\lam) v(R)r(R) \le n \le c_7(\lam) v(R)r(R).
\end{align}
\end{prop}

\proof
(1) (\ref{e:24-11}) is immediate by Proposition \ref{bkb}(a), and  
(\ref{e:24-12}) follows from Proposition \ref{nash}(a).\\
Using  Proposition~\ref{nash}(b), and writing
$f_n(y)=p_n(0,y)+p_{n+1}(0,y)$, $n'= 2 \lfloor n/2 \rfloor$,
\eq
 f_{n}(y) \le f_{n}(0) + |f_{n}(y)-f_{n}(0)| \le
 f_{n}(0)+ (c \Reff(0,y) n^{-1} p_{n'}(0,0) )^{1/2}.
\en
So, by (\ref{e:24-12}) and by the definition of $J(\lam)$, 
if $y\in B(R)$ then we have (\ref{e:24-13}), namely
\eq
 f_{n}(y) \le \frac{c'(\lam)}{v({\cal I}(n))}.
\en
\newline
(2)  Set $m=2\lam$, $\eps^{\al_1} =1/(2C_2m\lam)= 1/(4C_2\lam^2)$. 
Since $R \in J(\lam)$,
we have, for $x \in B(0,\eps R)$,
$$ \frac{r(R)}{\lam} \le \Reff(0, B^c) \le \Reff(0,x) + \Reff(x,B^c)
 \le \lam r(\eps R) + \Reff(x,B^c)\le \lam C_2\eps^{\al_1}r(R) + \Reff(x,B^c),$$
where we used (\ref{vrdoubl}) in the last inequality. 
Hence $\Reff(x,B^c) \ge r(R)/m$ if $x\in B(0,\eps R)$,
and so the assumption of Proposition \ref{bkb}(b) holds.
Since $R \in J(\lam)$, $V(R) \le \lam v(R)$.
Also $\frac12 \eps R= R/(2^{1+2/{\al_1}} C_2^{1/\al_1}\lam^{2/{\al_1}})=:c_2(\lam)R
\in J(\lam)$, so 
$V' \ge \lam^{-1} v(c_2(\lam)R)\ge c'(\lam)v(R)$ for some $c'(\lam)>0$; 
the bounds now follow
from  Proposition \ref{bkb}(b). 
\qed

\bigskip
Next we apply similar arguments to control $d(0,X_n)$, beginning with
a preliminary lemma.  Recall that $T_A$ was defined in \refeq{TAdef} to
be the hitting time of $A \subset G$.

\begin{lemma}\label{hitest} Let $\lam \ge 1$ and
$0< \eps^{\al_1} \le 1/(2C_2\lam^2)$. 
If $R \in J(\lambda)$, and $y \in B(\eps R)$ then
\begin{align}
P^y ( \tau_R < T_0) &\le  2 C_2\eps^{\al_1} \lam^2 ,\\
P^0 ( \tau_R < T_y) &\le   C_2\eps^{\al_1} \lam^2 .
\end{align}
\end{lemma}

\proof
If $A$ and $B$ are disjoint subsets of $G$ and $x \not\in A \cup B$, then
(see \cite[(4)]{BGP03})
\[ P^x(T_A < T_B) \le \frac{\Reff(x,B)}{\Reff(x,A)}. \]

Let $d(0,y) \le \eps R$. 
Then $\Reff(y,0)\le \lam r(d(y,0))\le \lam r(\eps R)\le \lam C_2\eps^{\al_1}r(R)$, 
while
\[ \Reff(y,B(R)^c) \ge \Reff(0,B(R)^c) - \Reff(0,y)
\ge {r(R)}/{\lam}-\lam C_2\eps^{\al_1}r(R) \ge r(R)/2\lam .\]
So,
\[ P^y(\tau_R < T_0) \le  \frac{\Reff(y,0)}{\Reff(y, B(R)^c)}
 \le 2C_2 \eps^{\al_1} \lam^2. \]
Similarly,
\[ P^0(\tau_R < T_y) \le  \frac{\Reff(0,y)}{\Reff(0, B(R)^c)}
 \le C_2\eps^{\al_1} \lam^2. \]
\qed

\begin{prop} \label{eqest2}
For each $\lambda>1$, there exist $c_1(\lam),\cdots,c_{10}(\lam)$
such that the following hold.\\
(a) Let 
$\eps < c_1(\lam)$ and $R, \eps R, c_2(\lam)\eps R
\in J(\lambda)$. Then
\eq
 \label{trdlb}
 P^y \big( \tau_R \le c_3(\lam) v(\eps R)r(\eps R) \big) \le c_4(\lam) \eps^{\al_1},
 \quad \text{ for } y \in B(\eps R).
\en
(b) Let $n\ge 1$, $M\ge 1$, and set $R=M {\cal I}(n)$.
If $R, c_5(\lam)R/M, c_6(\lam)R/M \in J(\lambda)$, then
\eq  P^0\big(\frac{d(0,X_n)}{{\cal I}(n)}>M\big)\le \frac{c_7(\lam)} {M^{\al_1}}. \en
(c) Let $R={\cal I}(n)$ and $\theta\in (0,1]$. If $R,\theta R \in J(\lam)$ then
\eq P^0\big( X_n \in B(\theta R) \big) \le c_8(\lam)\theta^{d_1}. \en
(d) Let $R={\cal I}(n)$.
If $R, c_9(\lam) R \in J(\lam)$ then
\eq P^0( \tau_{c_9(\lam) R} \le n) \ge
    P^0\big( X_n \not\in B(0,c_9(\lam)R)\big) \ge \tfrac12. \en
Hence
\eq
\lbeq{Edlb}
       E^0 d(0,X_n) \ge c_{10}(\lam) {\cal I}(n).
\en
\end{prop}

\proof
(a) 
Let $c_1(\lam)=(2^{1+1/\al_1}C_2^{1/\al_1}\lam^{2/\al_1})^{-1}\wedge 1$,
$c_2(\lam)=c_{\ref{eqest}.2}(\lam)$, and 
$c_3(\lam)=c_{\ref{eqest}.4}(\lam)< 1$. Then the desired
inequality is trivial when $\eps R\le 1$, so assume that $\eps R>1$.
Let  $q(y)=P^y(\tau_R < T_0)$, so that,
by substituting $2\eps$ into $\eps$ in Lemma \ref{hitest}, 
if $d(0,y) \le 2\eps R$ then $q(y)\le c_0 \eps^{\al_1} {\lam}^2$. 
Write $t_0=  c_3(\lam) v(\eps R)r(\eps R)$ and
$a=P^0(\tau_R\le t_0)$.  Now if $y \in B(2\eps R)$ then
\begin{align}
 P^y(  \tau_R \le t_0)
 &= P^y( \tau_R \le t_0, \tau_R < T_0) + P^y( \tau_R \le t_0,  \tau_R > T_0)\nnb
 &\le  P^y( \tau_R \le T_0) +  P^y(T_0< \tau_R,  \tau_R -T_0 \le t_0)\nnb
\label{pytub}
 &\le q(y) + (1-q(y)) a  \le c_0 \eps^{\al_1} \lam^2 + a,
\end{align}
using the strong Markov property for the second inequality.
So, by  a second application of
the strong Markov property, and (\ref{e:24-22}),
\begin{align}
a = P^0(\tau_R\le t_0) &\le
E^0 [1_{\{\tau_{\eps R}\le t_0 \}} P^{X_{\tau_{\eps R}}}( \tau_R \le t_0)] \nnb
  &\le  (1-c_{\ref{eqest}.5}(\lam) )
  (c_0 \eps^{\al_1} \lambda^2 +a),
\end{align}
where we used the fact
that $X_{\tau_{\eps R}}\in B(\eps R+1)\subset B(2\eps R)$
in the last inequality.
Rewriting this gives
$a \le c_0 \eps^{\al_1}\lam^2 (1-c_{\ref{eqest}.5}(\lam))/c_{\ref{eqest}.5}(\lam)$.
Substituting in (\ref{pytub}) gives (\ref{trdlb}) with
$c_4(\lam)=c_0\lam^2 /c_{\ref{eqest}.5}(\lam)$.
\\
(b)

Let $c_5(\lam)=c_*c_3(\lam)^{-1/(d_1+{\al_1})}$, 
where $c_*>0$ large 
is chosen later. Let  
$c_6(\lam)= c_2(\lam)c_5(\lam)$, $c_7(\lam)=(c_5(\lam)/c_1(\lam))^{\al_1}(c_4(\lam)\vee 1)$, $M'=M/c_5(\lam)$, and
$\eps=(M')^{-1}$. The desired
inequality is trivial when $c_7(\lam)/M^{\al_1}\ge 1$, so assume that
$c_7(\lam)/M^{\al_1}<1$. Then, $M>c_5(\lam)/c_1(\lam)$, so  
$\eps=c_5(\lam)/M<c_1(\lam)$. Thus the assumption in (a) is satisfied. 
Using (\ref{vrdoubl}), we have 
${\cal I}(n/c_3(\lam))\le 
{\hat c}c_3(\lam)^{-1/(d_1+\al_1)}
{\cal I}(n)$, so taking 
$c_*={\hat c}$, we have ${\cal I}(n/c_3(\lam))
\le \eps R$, which is equivalent to
\eq\label{add5-9}
n\le c_3(\lam)v(\eps R)r(\eps r). 
\en 
Since
\eqalign
 P^0(d(0,X_n)/{\cal I}(n)> M) & =
 P^0(d(0,X_n)  > R)\nnb
 &\le P^0(\tau_{R}\le n)
 \le P^0(\tau_R \le c_3(\lam) v(\eps R)r(\eps R)),
\enalign
where (\ref{add5-9}) is used in the last inequality. 
Using (a) gives the desired estimate.\\
(c)  By (\ref{e:24-13}),  writing $B'=B(0,\theta R)\subset B(0,R)$
and $f_n(0,y)=p_n(0,y)+p_{n+1}(0,y)$,
\eq\label{5-92tus}
  P^0( X_n\in B')= \sum_{y\in B'} p_n(0,y)\mu_y \le
   \sum_{y\in B'} f_n(0,y) \mu_y \le V(\theta R) c_{\ref{eqest}.1}(\lam)/ v(R).
\en
Since $\theta R\in J(\lam)$, using (\ref{vrdoubl}),
the right hand side of (\ref{5-92tus}) is bounded from above by 
$c_8(\lam)\theta^{d_1}$. \\
(d) Let $\theta=c_9(\lam)\in (0,1]$ satisfy
$c_8(\lam)\theta^{d_1} =\frac12$. Then, since $R,\theta R \in J(\lam)$, 
applying (c),
\eq P^0( X_n\in B(\theta R)) \le c_8(\lam)\theta^{d_1} =\tfrac12.  \en
This proves the first assertion. Also,
\eq
  E^0 d(0,X_n)\ge \theta R P^0( X_n \not\in B') \ge \tfrac12 \theta R
 \ge c_{10}(\lam) {\cal I}(n). \en\qed

\subsection{Proof of Propositions~\ref{ptight}--\ref{pmeans} and
Theorem~\ref{thm-rwre}}
\label{sec-genrwre}
We now consider a family of random graphs, as described in Section
\ref{sec-rwre}, and prove Propositions~\ref{ptight}--\ref{pmeans} and
Theorem~\ref{thm-rwre}.

We begin by obtaining tightness of 
 $E^0  \tau_R/(v(R)r(R))$, $v({\cal I}(n)) p_{2n}(0,0)$, and 
$d(0,X_n)/{\cal I}(n)$. In the following, we set 
$l(\lam)=c_{\ref{eqest}.2}(\lam)$. 

\medskip\noindent
\emph{Proof of Proposition \ref{ptight}}.
We begin with (\ref{pt-a}).
Let $\eps>0$. Choose $\lam\ge 1$ such that $2p(\lam)< \eps$ --
here $p(\lam)$ is the function given by Assumption~\ref{ass-rwre}.
Let $R\ge 1$ and set $F_1 =\{ R, l(\lam)R \in J(\lam)\}$.

Suppose first that $l(\lam)R \ge 1$. Then, by
Assumption~\ref{ass-rwre}(1), $\bP(F_1)\ge 1-2p(\lam)$.
For $\omega \in F_1$, by Proposition~\ref{eqest}, there exists $c_1<\infty$,
$q_1 \ge 0$ such that
\eq
 \label{e:ttz}
 (c_1 \lam^{q_1})^{-1} \le  E^x_\omega \tau_R/(v(R)r(R)) \le c_1\lam^{q_1}
 \text{ for } x \in B(l(\lam)R).
\en
So, if $\theta_0=  c_1 \lam^{q_1}$ then for $\theta \ge \theta_0$,
\eq
\label{e:tta}
\bP\big(  \theta^{-1} \le   E^0_\omega \tau_R/(v(R)r(R)) \le \theta  \big)
\ge \bP(F_1) \ge 1 -2 p(\lam) \ge 1 - \eps.
\en

Now consider the case when $R \le 1/l(\lam)$. For each graph
$\Gam(\omega)$ let
\[ Y(\omega) = \sup_{1 \le s \le 1/l(\lam)} E^0_\omega \tau_s/(v(s)r(s)). \]
Then $Y(\omega)<\infty$ for each $\omega$, so there
exists $\theta_1$ such that
\[ \bP(  E^0_\omega \tau_R/(v(R)r(R)) > \theta_1) \le \bP(Y > \theta_1)
 \le \eps. \]
If we take $\theta_1> v(1/l(\lam))r(1/l(\lam))$ then since
$E^0_\omega \tau_R \ge 1$, we have
$E^0_\omega \tau_R/(v(R)r(R))\ge \theta_1^{-1}$.
So, for $\theta\ge \theta_1$, we also have
$\bP\big(  \theta^{-1} \le  E^0_\omega \tau_R/(v(R)r(R))\le \theta  \big)
\ge 1 - \eps$, which completes the proof of  (\ref{pt-a}).

We now turn to (\ref{pt-b}).
Let $n \ge 1$, $\lam \ge 1$, and let $R_0$, $R_1$ be defined by
$n=c_{\ref{eqest}.7}(\lam) v(R_1)r(R_1)=v(R_0)r(R_0)$.
Let $F_2 =\{ R_0, R_1, l(\lam)R_1 \in J(\lam)\}$.
Suppose first that $R_0$ and $l(\lam)R_1$
are both greater than $1$; then $\bP(F_2) \ge 1-3p(\lam)$.
If $\omega \in F_2$ then by Proposition~\ref{eqest}
\[ (c_2 \lam^{q_2})^{-1} \le v({\cal I}(n)) p_{2n}^\omega(0,0) \le c_2\lam^{q_2},\]
for some $c_2>0, q_2>0$. So, 
\eq
\lbeq{ttaaa}
     \bP\big(  (c_2 \lam^{q_2})^{-1} \le v({\cal I}(n))p_{2n}^\omega
     (0,0) \le c_2 \lam^{q_2} \big)
 \ge \bP(F_2) \ge 1-3p(\lam).
\en
The case when $n$ is small is dealt with in the same way as in the proof
of (\ref{pt-a}).

Next we prove (\ref{pt-c}). Let $n \ge 1$ and $\lam\ge 1$.
Let $M=(\lam c_{\ref{eqest2}.7}(\lam))^{1/\al_1}=:l_1(\lam)$, and set
\eq
 \label{rdefa}
 R_0= M {\cal I}(n), \quad R_1= c_{\ref{eqest2}.5}(\lam){\cal I}(n), \quad
 R_2= c_{\ref{eqest2}.6}(\lam){\cal I}(n),
\en
$F_3 =\{ R_0, R_1, R_2 \in J(\lam)\}$.
Suppose first that $n$ is large enough so that
$R_i\ge 1$ for $0\le i \le 2$.
By Proposition~\ref{eqest2}(b), if $\omega \in F_3$ then
\[
    P^0_\omega\big(\frac{d(0,X_n)}{{\cal I}(n)}> l_1(\lam) \big)\le
    \frac{c_{\ref{eqest2}.7(\lam)}}{l_1(\lam)^{\al_1}}
    = \frac 1{\lam}.
\]
Taking $\theta=l_1(\lam)$, we have
\begin{eqnarray}
 P^* \big( \frac{d(0,X_n)}{{\cal I}(n)}> \theta \big) &\le & \bP(F_3^c) +
 \bE \big(  P^0_\omega\big(\frac{d(0,X_n)}{{\cal I}(n)}> l_1(\lam) \big);F_3\big)\nonumber\\
&\le &3p(l_1^{-1}(\theta)) + \frac 1{l_1^{-1}(\theta)},\label{e:dtta}
\end{eqnarray}
where $l_1^{-1}(\theta)$ is the inverse of $l_1(\theta)$. 

Now let $\eps>0$. Choose $\theta_0$ so that the right side of
(\ref{e:dtta}) is less than $\eps$. Let $l_1(\lam)=\theta_0$.
Then there exists $n_1=n_1(\eps)$ such that if $n\ge n_1$ then
$R_0, R_1, R_2$ (given by (\ref{rdefa})) are all greater than 1.
If $n\ge n_1$ then (\ref{e:dtta}) implies that
$P^*(d(0,X_n)/{\cal I}(n)> \theta_0 )<\eps$.

To handle the case when $n \le n_1$, for each $\omega$ let
\[ Z_\theta(\omega)
 = \max_{1\le n \le n_1} P^0_\omega(d(0,X_n)/{\cal I}(n)> \theta ). \]
Then $Z_\theta$ is non-increasing in $\theta$, and
$\lim_{\theta \to \infty} Z_\theta(\omega) =0$ for each $\omega$.
So, by monotone convergence
\[ \lim_{\theta \to \infty} \bE  Z_\theta(\omega) =0. \]
Thus there exists $\theta_1$ such that
\[ P^*(d(0,X_n)/{\cal I}(n)> \theta_1 )\le \bE Z_{\theta_1} < \eps
\quad \text { for all }  n\le n_1. \]
Taking $\theta = \theta_0 \vee \theta_1$,
we obtain  (\ref{pt-c}).

Finally, we prove (\ref{pt-d}).
Let $\eps>0$. 
Choose $\lam$ so that $2p(\lam)+ 1/c_{\ref{eqest2}.8}(\lam)< \eps$, and let
$\theta_0= c_{\ref{eqest2}.8}(\lam)^{2/d_1}$, $\delta=1/\theta_0$.
Choose $R$ so that  $v(R)r(R)=n$, and $n_0=n_0(\eps)$ such that
$n\ge n_0$ implies $\delta R \ge 1$. Set
$\theta_1= 1+{\cal I}(n_0)$, and $\theta = \theta_0 \vee \theta_1$.
Suppose first that $n \ge n_0$, and set
$F_4 = \{ R, \delta R \in J(\lam)\}$. If $\omega \in F_4$
then by Proposition~\ref{eqest2}(c), we have
\[
 P^0_\omega\big(d(0,X_n)/{\cal I}(n) \le \delta \big)
\le c_{\ref{eqest2}.8}(\lam) \delta^{d_1}. \]
So,
\begin{align}
 P^* \big( (1+ d(0,X_n))/{\cal I}(n)<  \theta^{-1} \big) &\le
P^*\big(d(0,X_n)/{\cal I}(n)< \theta_0^{-1} \big) \nnb
 &\le \bP(F_4^c)
   +\bE \big(  P^0_\omega( d(0,X_n)/{\cal I}(n)<   \theta_0^{-1} );F_4\big)\nnb
 &\le 2p(\lam) + 1/c_{\ref{eqest2}.8}(\lam)\le \eps.
\label{e:dttb}
\end{align}
If $n \le n_0$ then
$(1+d(0,X_n))/{\cal I}(n) \ge 1/{\cal I}(n) \ge \theta_1^{-1}$, and so
we deduce that, for all $n$,
\[ P^*((1+ d(0,X_n))/{\cal I}(n)<  \theta^{-1}) \le \eps, \]
which proves (\ref{pt-d}). \qed

\medskip\noindent
\emph{Proof of Proposition~\ref{pmeans}.}
We begin with the upper bounds in 
(\ref{e:mmean}). 
By (\ref{e:newads}) and Assumption~\ref{ass-rwre}(2),
\[ \bE (E^0_\omega \tau_R) \le \bE (\Reff(0,B(R)^c)V(R)) \le c v(R)r(R).\]

For the lower bounds, it is sufficient to find a set $F\subset \Omega$
of `good' graphs with $\bP(F)\ge c>0$ such that, for all $\omega \in F$
we have suitable lower bounds on $E^0_\omega \tau_R$, $p^\omega_{2n}(0,0)$
or $E^0_\omega d(0,X_n)$.
For the lower bounds,
we assume that $R\geq 1$ is large enough so that $l(\lam_0)R\ge 1$,
where $\lam_0$ is
chosen large enough that $p(\lambda_0)< 1/8$.
We can then obtain the results for all
$n$ (chosen below to depend on $R$) and $R$ by
adjusting the constants $c_1,\cdots, c_4$ in (\ref{e:mmean})--(\ref{e:dmean}).

Let
$F =\{ R, l(\lam_0)R \in J(\lam_0)\}$. Then $\bP(F)\ge \frac34$,
and for $\omega\in F$, by (\ref{e:24-21}), $E^0_\omega \tau_R
\ge c_1(\lam_0) v(R)r(R)$. So,
\[ \bE(E^0_\omega\tau_R)\ge \bE(E^0_\omega\tau_R;F)
\ge c_1(\lam_0) v(R)r(R) \bP(F) \ge c_2(\lam_0) v(R)r(R). \]
Also, by (\ref{e:24-23}), if
$n\in [\frac12 c_{\ref{eqest}.7}(\lam_0) v(R)r(R), c_{\ref{eqest}.7}(\lam_0) v(R)r(R)]$
then
\[ p_{2n}^\omega(0,0) \ge \frac{c_3(\lam_0)}{v({\cal I}(n))}. \]
Given $n \in \bN$, choose $R$ so that $n=c_{\ref{eqest}.7}(\lam_0) v(R)r(R)$ and let
$F$ be as above. Then
\[ \bE p_{2n}^\omega(0,0) \ge \bP(F)\frac{c_3(\lam_0)}{v({\cal I}(n))} \ge 
\frac{c_4(\lam_0)}{v({\cal I}(n))}, \]
giving the lower bound in (\ref{e:pmean}).

A similar argument uses \refeq{Edlb} to conclude (\ref{e:dmean}). 

Finally we prove (\ref{e:pmeanII}). 
Let $H_k$ be the event of the left hand side of (\ref{newcond45}) with 
$\lam=k$. 
By Proposition \ref{nash}(a), we see that 
$p_{2n}^\omega(0,0)\le c_1k/v({\cal I}(n))$ if 
$\omega \in H_k$,
where $R$ is chosen to satisfy $v(R)r(R)/2\le n\le v(R)r(R)$. 
Since
$\bP(\cup_k H_k)=1$, using 
(\ref{newcond45}), we have 
\begin{eqnarray*}
\bE p_{2n}^\omega(0,0) &\le&\sum_k\frac{c_1(k+1)}{v({\cal I}(n))}
\bP(H_{k+1}\setminus H_k)
\le \sum_k\frac{c_1(k+1)}{v({\cal I}(n))}
\bP(H_k^c)\\
&\le & \frac{c_2}{v({\cal I}(n))}\sum_k(k+1)k^{-q_0'}<\infty,
\end{eqnarray*}
since $q_0'>2$. We thus obtain (\ref{e:pmeanII}). \qed

\medskip\noindent
\emph{Proof of Remark~\ref{resremvi}.\,1.}
In this case, we have 
\[\bP(\{\Reff (0,y)\le  r(d(0,y)),~\forall y\in B(R)\})=1,\]
so, similarly to the proof of Proposition \ref{nash}, 
for $v(R)r(R)/2\le n\le v(R)r(R)$, we have 
$$ f_{2n}(0) \le c_1(\frac {1}{V(0,R)}\vee \frac 1{v(R)})
\le \frac{c_2}{v({\cal I}(n))} ( \frac{v(R)}{V(0,R)} \vee 1). $$
Using this and (\ref{e:addinvvol}), we have 
\[ \bE p_{2n}^\omega(0,0) 
 \le \frac c{v({\cal I}(n))} \bE (1 + \frac{v(R)}{V(0,R)} )\le 
 \frac{c'}{v({\cal I}(n))},\]
so 
(\ref{e:pmeanII}) is obtained. \qed

\medskip\noindent
\emph{Proof of Theorem~\ref{thm-rwre}.}
(I) We will take $\Omega_0 = \Omega_a \cap \Omega_b \cap \Omega_c$
where the sets $\Omega_*$ are defined in the proofs of (a), (b) and (c).
By  Assumption~\ref{ass-rwre}(3),
$p(\lam) = \bP( R \not\in J(\lam))\le c_0 \lam^{-q_0}$.

\noindent (a) We begin with the case $x=0$, and
write $w(n) = p^\omega_{2n}(0,0)$. By \refeq{ttaaa} we have
\[ \bP( (c_1 \lam^{q_1})^{-1} < v({\cal I}(n)) w_n \le c_1 \lam^{q_1})
 \ge 1-3 p(\lam). \]
Let $n_k= \lfloor e^k \rfloor $ and $\lam_k= k^{2/q_0}$. Then,
since $\sum p(\lam_k) < \infty$, by Borel--Cantelli
there exists $K_0(\omega)$ with $\bP(K_0<\infty)=1$ such that
$c_1^{-1} k^{- 2 q_1/q_0} \le v({\cal I}(n_k)) w(n_k)
 \le c_1  k^{2 q_1/q_0}$
for all $k\ge K_0(\omega)$. Let $\Omega_a =\{ K_0 < \infty\}$.
For $k \ge K_0$ we therefore have
\[  c_2^{-1} \frac{(\log n_k)^{- 2 q_1/q_0}}{v({\cal I}(n_k))} \le w(n_k)
 \le c_2 \frac{(\log n_ k)^{2 q_1/q_0}}{v({\cal I}(n_k))}, \]
so that (\ref{e:logpnlima}) holds for the subsequence $n_k$.
The spectral decomposition gives that $p^\omega_{2n}(0,0)$ is monotone
decreasing in $n$.
So, if $n >N_0= e^{K_0} +1$, let $k \ge K_0$ be such that
$n_k \le n < n_{k+1}$.
 Then
\[ w(n) \le w(n_k) \le c_2 \frac{(\log n_k)^{2 q_1/q_0}}{v({\cal I}(n_k))}
 \le c_2'  \frac{(\log n)^{2 q_1/q_0}}{v({\cal I}(n))}. \]
Similarly $w(n) \ge w(n_{k+1}) \ge \frac{c_3}{v({\cal I}(n))} (\log n)^{-2q_1/q_0}$.
Taking $q_2 > 2q_1/q_0$, so that the constants $c_2, c_3$ can be
absorbed into the $\log n$ term, we obtain
\eq
 \label{p2nlim}
  \frac{(\log n)^{-q_2}}{v({\cal I}(n))} \le  p^\omega_{2n}(0,0)
 \le  \frac{(\log n)^{q_2}}{v({\cal I}(n))} \quad \text { for all }  n\ge N_0(\omega).
\en

If $x, y \in \Ccal(\omega)$ and $k=d_\omega(x,y)$, then
using the Chapman--Kolmogorov equation
\[ p^\omega_{2n}(x,x)(p^\omega_{k}(x,y)
\mu_x(\omega))^2 \le p^\omega_{2n+2k}(y,y). \]
Let $\omega \in \Omega_a$, $x \in \Ccal(\omega)$, write $k=d_\omega(0,x)$,
$h^\omega(0,x)=(p^\omega_{k}(x,0)\mu_x(\omega))^{-2}$,
and let $n \ge N_0(\omega) + 2k$.
Then
\begin{align*}
 p^\omega_{2n}(x,x) &\le h^\omega(0,x) p^\omega_{2n+2k}(0,0) \\
 &\le  h^\omega(0,x)  \frac{(\log(n+k))^{q_2}}{v({\cal I}(n+k))}\\
 &\le  h^\omega(0,x)  \frac{(\log (2n))^{q_2}}{v({\cal I}(n))}
 \le  \frac{(\log n)^{1+q_2}}{v({\cal I}(n))} 
\end{align*}
provided $ \log n \ge 2^{q_2} h^\omega(0,x)$.
Taking 
\eq\label{e:nxallx}
N_x(\omega) = \exp (2^{q_2} h^\omega(0,x)) + 2d_\omega(0,x) + N_0(\omega),
\en
and $\beta_1=1+q_2$, this gives the upper bound in (\ref{e:logpnlima}).
The lower bound is obtained in the same way.

\medbreak
\noindent (b) Let $R_n = e^n$ and $\lam_n = n^{2/q_0}$.
Let $F_n =\{ R_n, l(\lam_n)R_n \in J(\lam_n) \}$.
Then (provided $l(\lam_n)R_n \ge 1$) we have
$\bP(F_n^c) \le 2 p(\lam_n) \le 2 n^{-2}$.
So, by Borel--Cantelli, if $\Omega_b =\liminf F_n$, then
$\bP(\Omega_b)=1$. Hence there exists $M_0$ with $M_0(\omega)<\infty$
on $\Omega_b$, and such that $\omega \in F_n$ for all $n \ge M_0(\omega)$.

Now fix $\omega \in \Omega_b$, and let $x \in \Ccal(\omega)$.
Write $F(R) = E^x_\omega \tau_R$. By (\ref{e:ttz}) there exist
constants $c_4$, $q_4$ such that
\eq \label{mrnbound}
 (c_4 \lam_n^{q_4})^{-1} \le \frac{F(R_n)}{v(R_n)r(R_n)} \le c_4 \lam_n^{q_4} .
\en
provided $n \ge M_0(\omega)$ and $n$ is also large enough so that 
$x \in B(l(\lam_n)R_n)$. Writing $M_x(\omega)$ for the smallest such $n$,
\[ c_4^{-1} (\log R_n)^{-2 q_4/q_0} v(R_n)r(R_n)
  \le F(R_n) \le  c_4 (\log R_n)^{2 q_4/q_0} v(R_n)r(R_n),
\quad \text{ for all } n \ge M_x(\omega). \]
As $F(R)$ is monotone increasing, the same argument as in (a) enables
us to replace $F(R_n)$ by $F(R)$, for all $R\ge R_x= 1+ e^{M_x}$.
Taking $\beta_2 > 2q_4/q_0$ we obtain (\ref{e:logtaulima}).

\medbreak
\noindent (c) Recall that $Y_n =\max_{0\le k \le n} d(0,X_k)$.
We begin by noting that
\eq
 \label{ytaurel}
   \{ Y_n \geq R \} =\{ \tau_R \le n\}.
\en
Using this, (\ref{e:ynlim}) follows easily from (\ref{e:rnlim}).

It remains to prove (\ref{e:rnlim}). Since $\tau_R$ is monotone
in $R$, as in (b) it is enough to prove the result
for the subsequence $R_n=e^n$.

The estimates in (b) give the upper bound. In
fact, if $\omega \in \Omega_b$, and
$n \ge M_x(\omega)$, then
by (\ref{mrnbound})
\[ P^x_\omega( \tau_{R_n} \ge n^2 c_4 \lam_n^{q_4}v(R_n)r(R_n)) \le
  \frac{F(R_n)}{ n^2 c_4 \lam_n^{q_4} v(R_n)r(R_n)} \leq n^{-2}. \]
So, by Borel--Cantelli (with respect to the law $P^x_\omega$),
there exists $N'_x(\omega, \overline \omega)$ with
\[ P^x_\omega(N'_x< \infty)
 =P^x_\omega( \{ \overline \omega: N'_x(\omega, \overline \omega)<\infty\})=1 \]
such that
\[ \tau_{R_n} \le c_5 (\log R_n)^{q_5} v(R_n)r(R_n), \quad
\text{ for all } n \ge N'_x. \]

For the lower bound, write $c_{\ref{eqest2}.3}(\lam)=c_6 \lam^{-q_6}$,
$c_{\ref{eqest2}.4}(\lam)=c_7 \lam^{q_7}$,  
where we choose $q_6+q_7\ge 2$.  
Let $\lam_n= n^{2/q_0}$, and 
$\eps_n^{\al_1}= cn^{-2} \lam_n^{-q_6 - q_7}$. Here $c>0$ is chosen 
small enough so that $\eps_n\le c_{\ref{eqest2}.1}(\lam)$. 
Set $G_n=\{ R_n, \eps_n R_n, l(\lam_n)\eps_nR_n \in J(\lam_n)\}$.
Then, for $n$ sufficiently large so that $l(\lam_n)\eps_nR_n \ge 1$,
we have $\bP(G_n^c) \le 3 p(\lam_n) \le 3 c_0 n^{-2}$.
Let $\Omega_c =\Omega_b \cap (\liminf G_n)$; then by Borel--Cantelli, 
$\bP(\Omega_c)=1$ and there exists $M_1$ with $M_1(\omega)<\infty$
for $\omega \in \Omega_c$ such that $\omega \in G_n$ whenever
$n \ge M_1(\omega)$.
By Proposition \ref{eqest2}(a), if $n\ge M_1$ and
$x \in B(\eps_n R_n)$ then
\eq
 P^x_\omega( \tau_{R_n} \le c_6 \lam_n^{-q_6} v(\eps_n R_n)r(\eps_n R_n))
 \le c_7 \lam_n^{q_7} \eps_n^{\al_1} \le c_7' n^{-2}.
\en
So, using Borel--Cantelli, we deduce that (for some $q_8$)
\[ \tau_{R_n} \ge  c_6 \lam_n^{-q_6} v(\eps_n R_n)r(\eps_n R_n)
 \ge n^{-q_8} v(R_n)r(R_n) = (\log R_n)^{-q_8} v(R_n)r(R_n), \]
for all $n \ge N''_x(\omega, \overline \omega)$.
This completes the proof of (\ref{e:rnlim}).

\medbreak

The proof of (II) is similar by the following changes; take $\lam_k=(e+(2/c_4) 
\log k)^{1/q_0}$ instead of $\lam_k=k^{2/q_0}$, and take 
$N_x(\omega) = \exp(\exp (C h^\omega(0,x))) + 2d_\omega(0,x) + N_0(\omega)$
in (\ref{e:nxallx}). Then, $\log n$
(resp. $\log n_k, \log R_n$) in the above proof are changed to $\log\log n$
(resp. $\log\log n_k, \log\log R_n$) and the proof of (a) and (b) goes through. Since 
the modifications are simple, we omit details. 

\medbreak

We now prove (III). 
For (a), $\lim_n \log p_{2n}^\omega(0,0)/\log n = -D/(D+{\al})$, $\bP$--a.s. is 
easy from (\ref{e:logpnlima}) and (\ref{logvrcond}). 
Since $\sum_n p^\omega_{2n}(0,0)=\infty$, $X$ is recurrent.
(b) is also easy from (\ref{e:logtaulima}) and (\ref{logvrcond}). 

\noindent
(c) We first consider the case $x=0$.
Let $c_1\in (0,1)$, $c_2\ge 2$, $q_1\ge 1$, 
$q_2\ge 2$ be chosen so that
\[   c_{\ref{eqest2}.3}(\lam) \ge c_1 \lam^{-q_1}, \quad
     c_{\ref{eqest2}.4}(\lam) \le c_2 \lam^{q_2}.
\]

Let $R_k=e^k$, and $\lam_k=k^{q_3}$ where $q_3\ge 2$ is chosen large
enough so that $\sum p(\lam_k) < \infty$. Let
$\eps_k^{\al_1}= c_2^{-1} \lam_k^{-q_2} k^{-q_3}$.
Set
\[ F_k = \{ R_k, \eps_k R_k, c_{\ref{eqest2}.2}(\lam)\eps_k R_k\in J(\lam_k) \}. \]
For $\omega \in F_k$ we have
by Proposition \ref{eqest2}(a)
\[ P^0_\omega( \tau_{R_k} \le c_1 \lam_k^{-q_1} v(\eps_k R_k)r(\eps_k R_k))
  \le c_2 \lam_k^{q_2} \eps_k^{\al_1} = k^{-q_3}. \]
Set $n(k)= c_1 \lam_k^{-q_1} v(\eps_k R_k)r(\eps_k R_k)\ge c_3\lam_k^{-q_1}
(\eps_k R_k)^{D+\al}(\log (\eps_k R_k))^{-m_1-m_2}$. Then
\eq
  P^*(  \{\tau_{R_k} \le n(k)\} \cup F_k^c ) \le \bP(F_k^c) +  k^{-q_3}
\le 3p(\lam_k) +   k^{-q_3}.
\en
Therefore by Borel--Cantelli, we deduce that, $P^*$--a.s.,
for all sufficiently large $k$, $\tau_{R_k} > n(k)$ and $F_k$
holds. So, for large $k$,
\[   S_{n(k)} \le S_{\tau_{R_k}} \le V(R_k) \le \lam_k v(R_k)
\le c_4\lam_kR_k^D(\log R_k)^{m_1}. \]
If $n$ is sufficiently large, then choosing $k$ so that
$n(k-1) < n \le n(k)$,
\begin{eqnarray*}
  \frac{\log S_n}{\log n} &\le&
\frac{\log S_{n(k)}}{\log n(k-1)}
 \le \frac{Dk + \log (c_4\lam_k)+m_1\log k} {(D+{\al})(k-1)+\log(c_3 
 \eps_{k-1}^{D+{\al}}  \lam_{k-1}^{-q_1})-(m_1+m_2)\log\log(\eps_ke^k)}\\
& \le& \frac k{k-1}\frac D{D+{\al}} + \frac{c_5\log k}{k}, 
\end{eqnarray*}
and this gives the upper bound in (\ref{e:snlim}) for the case $x=0$.

For the lower bound, let
$\xi(x,R) = 1_{\{ T_x > \tau_R\}}$.
If $R \in J(\lam)$ and $\eps^{\al_1}< 1/(2C_2\lam^2)$ then by Lemma \ref{hitest},
\[ P^0_\omega( \xi(x,R) =1 ) \le C_2\eps^{\al_1} \lam^2, \quad
\text{ for } x \in B(\eps R). \]
Set
\[ Y_k = V(\eps_k R_k)^{-1} \sum_{x \in B(\eps_k R_k)} \xi(x,R_k)\mu_x. \]
Then if $\omega \in F_k$,
\[ P^0_\omega( Y_k \ge \tfrac12) \le 2 E^0_\omega Y_k
\le 2 C_2\eps_k^{\al_1} \lam_k^2 \le c_6k^{-q_3}. \]
Let $m(k) =  k^{q_3} \lam_k^2 v(R_k)r(R_k)\le c_7
k^{q_3} \lam_k^2 R_k^{D+\al}(\log R_k)^{m_1+m_2}$. Then  if $\omega \in F_k$,
by (\ref{e:24-11}),
\[  P^0_\omega( \tau_{R_k} \ge m(k) ) \le 2 \lam_k^2 v(R_k)r(R_k) m(k)^{-1}
 = 2 k^{-q_3}. \]
Thus
\[ P^*( F_k^c \cup \{ Y_k \ge \tfrac12\} \cup\{  \tau_{R_k} \ge m(k)\} )
\le 3p(\lam_k) + (2+c_6) k^{-q_3}, \]
so by Borel--Cantelli, $P^*$--a.s. there exists a $k_0(\omega)<\infty$
such that, for all $k\ge k_0$, $F_k$ holds, $\tau_{R_k} \le m(k)$,
and $Y_k \le 1/2$. So, for $k\ge k_0$,
\begin{eqnarray*} S_{m(k)} &\ge&
S_{\tau_{R_k}} = \sum_{x \in B(\eps_k R_k)} (1-\xi(x,R_k))\mu_x=
  V(\eps_k R_k) (1-Y_k)\\ 
  &\ge &\tfrac12 \lam_k^{-1} v(\eps_k R_k)\ge c_8 \lam_k^{-1}
  (\eps_k R_k)^D(\log (\eps_k R_k))^{-m_1}
  \end{eqnarray*}
Let $n$ be large enough so that $m(k) \le n < m(k+1)$ for some
$k \ge k_0$. Then
\[ \frac{\log S_n}{\log n} \ge \frac{ \log S_{m(k)}}{ \log m(k+1)}
  \ge \frac{ Dk - c \log k}{(D+{\al})(k+1) + c' \log (k+1)}, \]
and the lower bound  in (\ref{e:snlim}) follows.
This proves  (\ref{e:snlim}) when $x=0$.

Now let
\[ \Omega_0 = \{ \omega: G(\omega) \text{ is recurrent and }
 P^0_\omega( \lim_n (\log S_n/\log n) = \tfrac D{D+{\al}})=1 \}. \]
We have $\bP(\Omega_0)=1$. If $\omega \in \Omega_0$, and
$x \in G(\omega)$ then $X$ hits $0$ with $P^x_\omega$--probability
1. Since the limit does not depend on the initial segment
$X_0, \dots , X_{T_0}$, we obtain  (\ref{e:snlim}).

\qed


\begin{thebibliography}{10}
\bibitem{AN86}
Aizenman, M., and Newman, C. M. (1986).
\newblock Discontinuity of the percolation density in one dimensional 
$1/{|x-y|}^2$ Percolation Models.
\newblock {\em Commun. Math. Phys.} {\bf 107}, 611--647. 

\bibitem{AO82}
~Alexander, S., and ~Orbach, R. (1982).
\newblock Density of states on fractals: ``fractons''.
\newblock {\em J. Physique (Paris) Lett.} {\bf 43}, L625--L631. 

\bibitem{HBA00} 
Ben-Avraham, D., and Havlin, S. (2000). 
\newblock  Diffusion and Reactions in Fractals and Disordered Systems. 
\newblock Cambridge University Press, Cambridge. 

\bibitem{Barl04}
Barlow, M. T. (2004).
\newblock Random walks on supercritical percolation clusters.
\newblock {\em Ann. Probab.} {\bf 32}, 3024--3084.

\bibitem{BCK05}
Barlow, M. T., ~Coulhon, T., and ~Kumagai, T. (2005).
\newblock Characterization of sub-{Gaussian} heat kernel estimates on strongly
  recurrent graphs.
\newblock {\em Comm. Pure Appl. Math.} {\bf 58}, 1642--1677. 

\bibitem{BJKS}
Barlow, M. T., J\'arai, A. A., Kumagai, T., and Slade, G. (2008).
\newblock Random walk on the incipient infinite cluster for oriented 
percolation in high dimensions. 
\newblock {\em Commun. Math. Phys.} {\bf 278}, 385--431. 

\bibitem{BK06}
Barlow, M. T., and Kumagai, T. (2006).
\newblock Random walk on the incipient infinite cluster on trees. 
\newblock {\em Illinois J. Math.} {\bf 50}, 33--65 (Doob volume).

\bibitem{BBY07}
Benjamini, I., Berger, N., and Yadin, A. (2007).
\newblock Long range percolation mixing time.
\newblock {Preprint available at arXiv:math/0703872}. 

\bibitem{Ber02} 
~Berger, N. (2002).
\newblock Transience, recurrence and critical behavior
for long-range percolation.
\newblock {\em Commun. Math. Phys.} {\bf 226}, 531--558.

\bibitem{BB06}
~Berger, N., and ~Biskup, M. (2007).
\newblock Quenched invariance principle for simple random walk on percolation
  clusters.
\newblock {\em Prob. Theory Related Fields.} {\bf 137}, 83--120. 

\bibitem{BGP03}
~Berger, N., ~Gantert, N., and ~Peres, Y. (2003).
\newblock The speed of biased random walk on percolation clusters.
\newblock {\em Probab. Theory Related Fields.} {\bf 126}, 221--242.

\bibitem{CroKum}
~Croydon, D., and ~Kumagai, T. (2007). 
\newblock Random walks on Galton-Watson trees with infinite variance 
offspring distribution conditioned to survive.  
\newblock {Preprint available at http://www.math.kyoto-u.ac.jp/\~{ }kumagai/kumpre.html}.

\bibitem{Kest86}
~Kesten, H. (1986).
\newblock The incipient infinite cluster in two-dimensional percolation.
\newblock {\em Probab. Theory Related Fields.} {\bf 73}, 369--394. 

\bibitem{Kest86a}
~Kesten, H. (1986).
\newblock Subdiffusive behavior of random walk on a random cluster.
\newblock {\em Ann.\ Inst.\ H.\ Poincar\'e Probab.\ Statist.} {\bf 22}, 
425--487. 

\bibitem{Kiga01}
~Kigami, J. (2001).
\newblock {\em Analysis on Fractals}.
\newblock Cambridge University Press, Cambridge.

\bibitem{KozNach}
~Kozma, G., and ~Nachmias, A. (2008). 
\newblock The Alexander-Orbach conjecture holds in high dimensions. 
\newblock Preprint available at arXiv:0806.1442. 

\bibitem{MP06}
~Mathieu, P., and ~Piatnitski, A. (2007).
\newblock Quenched invariance principles for random walks on percolation
  clusters.
\newblock {\em Proceedings of the Royal Society A.} {\bf 463}, 2287--2307.  

\bibitem{Mis07}
Misumi, J. (2008). 
\newblock Estimates on the effective resistance in a long-range 
percolation on ${\bZ}^d$. 
\newblock  {\em J. Math. Kyoto Univ.}, to appear. 

\bibitem{SS04}
~Sidoravicius, V., and Sznitman, A. -S. (2004).
\newblock Quenched invariance principles for walks on clusters of percolation
  or among random conductances.
\newblock {\em Probab. Theory Related Fields}. {\bf 129}, 219--244. 

\end{thebibliography}
\end{document}